\documentclass[12pt]{article}
\usepackage{amssymb}
\usepackage{amscd}
\usepackage{amsmath}
\usepackage{amsthm}
\usepackage{verbatim}

\let\geq\geqslant
\let\leq\leqslant

\renewcommand{\a}{{\alpha}}

\renewcommand{\d}{{\delta}}
\newcommand{\e}{{\varepsilon}}

\newcommand{\s}{{\sigma}}

\DeclareMathOperator{\dist}{dist}
\DeclareMathOperator{\supp}{supp}

\DeclareMathOperator{\Ran}{Ran}
\DeclareMathOperator{\Dom}{Dom}
\DeclareMathOperator{\Ker}{Ker}

\newcommand{\abs}[1]{\lvert#1\rvert}

\newcommand{\norm}[1]{\lVert#1\rVert}

\newcommand{\R}{{\mathbb R}}
\newcommand{\N}{{\mathbb N}}

\newcommand{\C}{{\mathbb C}}

\renewcommand{\H}{{\mathcal H}}

\numberwithin{equation}{section}

\renewcommand{\P}{{\mathcal P}}
\theoremstyle{plain}
\newtheorem{theorem}{\bf Theorem}[section]
\newtheorem{lemma}[theorem]{\bf Lemma}

\theoremstyle{definition}

\theoremstyle{remark}
\newtheorem*{remark*}{\bf Remark}
\newtheorem{remark}[theorem]{\bf Remark}

\renewcommand{\qed}{\vrule height7pt width5pt depth0pt}

\newcommand{\Norm}[1]{{|\!|\!|}#1{|\!|\!|}}

\DeclareMathOperator{\Cp}{Cap}
\DeclareMathOperator{\codim}{codim}
\DeclareMathOperator{\esssup}{ess\, sup}
\DeclareMathOperator{\essinf}{ess\, inf}
\DeclareMathOperator{\osc}{osc}

\begin{document}

\date{3 April 2005}
\title{Spectral asymptotics of Pauli operators and 
orthogonal polynomials in complex domains}

%\date{}
\author{N.~Filonov\thanks{Department of Mathematical Physics, Faculty of Physics,
St.Petersburg State University, 198504 St.Petersburg, Russia.
email: filonov@mph.phys.spbu.ru}\,\,
and A.~Pushnitski\thanks{Department of Mathematical Sciences,
Loughborough University, Loughborough, LE11 3TU, U.K.
email: a.b.pushnitski@lboro.ac.uk}
}
\maketitle
\begin{abstract}
We consider the spectrum of a two-dimensional Pauli operator with 
a compactly supported electric potential
and a variable  magnetic field with a positive mean value.
The rate of accumulation of eigenvalues to zero
is described in terms of the logarithmic capacity
of the support of the electric potential. A connection between 
these eigenvalues and orthogonal polynomials in complex domains 
is established.
\end{abstract}

\emph{Keywords:} Pauli operator, magnetic field, spectral asymptotics, 
logarithmic capacity, orthogonal polynomials

%%%%%%%%%%%%%%%%%%%%%%%%%%%%%%%%%%%%%%%%
\section{Introduction}
%%%%%%%%%%%%%%%%%%%%%%%%%%%%%%%%%%%%%%%%
\textbf{1. The unperturbed Pauli operator.}
Let $B=B(x)$, $x=(x_1,x_2)\in\R^2$, be a 
real valued function which has the physical meaning of 
the strength of a magnetic field in $\R^2$.
A  two-dimensional non-relativistic spin-$1/2$ 
particle in the external magnetic field $B$
can be described by the Pauli operator
$$
h=
\begin{pmatrix}
h^+&0\\ 0 & h^-
\end{pmatrix}
\quad \text{in }L^2(\R^2)\oplus L^2(\R^2).
$$
The standard approach to the definition of the operators $h^\pm$ 
in $L^2(\R^2)$ involves
introducing the magnetic vector potential $A(x)=(A_1(x),A_2(x))$ 
such that $B=\partial_{x_1}A_2-\partial_{x_2}A_1$ 
and setting
\begin{equation}
h^\pm=(-i\nabla-A)^2\mp B.
\label{f8}
\end{equation}
Instead, we adopt the approach advocated in \cite{ErdosV},
which consists of defining $h^\pm$ in terms of a solution 
$\Psi=\Psi(x)$ to the differential equation $\Delta\Psi=B$.
\emph{Assume} that $B$ is such that a solution $\Psi$ 
can be chosen subject to the condition
\begin{equation}
\Psi(x)=\dfrac{B_0}{4}\abs{x}^2+\Psi_1(x), 
\quad \overline{\Psi}_1=\Psi_1\in L^\infty(\R^2),
\quad B_0>0.
\label{psi}
\end{equation}
Important examples of magnetic fields $B$ of this class are periodic 
fields with mean value $B_0$ and constant magnetic fields
$B(x)=B_0$.

Next, denote, as usual,
$\partial=\tfrac12(\partial_ {x_1}-i\partial_{x_2})$
and 
$\overline{\partial}=\tfrac12 (\partial_{x_1}+i\partial_{x_2})$.
Consider  the quadratic forms
\begin{equation}
{\mathfrak h}^+[u]=4\int_{\R^2}\abs{\overline{\partial}( e^{\Psi(x)} u(x))}^2 e^{-2\Psi(x)}dx,
\quad
{\mathfrak h}^-[u]=4\int_{\R^2}\abs{\partial (e^{-\Psi(x)} u(x))}^2 e^{2\Psi(x)}dx,
\label{f4}
\end{equation}
which are closed on the domains 
$\Dom ({\mathfrak h}^\pm)=\{ u\in L^2(\R^2)\mid {\mathfrak h}^\pm[u]<\infty\}$.
Let us define $h^\pm$ as the self-adjoint
operators in $L^2(\R^2)$, corresponding to the quadratic forms 
${\mathfrak h}^\pm$.
For a wide class of magnetic fields 
this definition is equivalent to the standard
definition \eqref{f8} with $A=(-\partial_{x_2}\Psi,\partial_{x_1}\Psi)$;
see \cite{ErdosV} for a detailed analysis of this issue.
In fact, the magnetic field $B$ or the magnetic vector 
potential $A$ do not enter directly either the 
definition of Pauli operator or any of our considerations; 
instead, the ``potential function" $\Psi$ becomes the main functional
parameter.
Note that the condition \eqref{psi} is very close to the `admissibility' 
condition used in \cite{Raikov3}.

We will denote by ${\mathfrak h}_0^\pm$ and $h_0^\pm$
the above defined forms and operators corresponding to 
the case of the constant magnetic field $B(x)=B_0>0$.

\textbf{2. Zero modes and the spectral gap. }
It is well known that Pauli operator $h$ has infinite dimensional kernel.
More precisely (cf \cite{AharonovC}) we have $\Ker h^-=\{0\}$ and 
\begin{equation}
\Ker h^+=\{u\in L^2(\R^2)\mid 
u(x)=f(x) e^{-\Psi(x)}, \overline\partial f=0\},
\quad \dim \Ker h^+=\infty.
\label{a16}
\end{equation}
Next, the following well known supersymmetric argument 
(which has appeared in many forms in the literature; see e.g. 
\cite{Iwatsuka} or \cite{Raikov3}) establishes the existence
of a spectral gap $(0,m)$, $m>0$ of the operator $h$.
Let $a_0$ and $a_0^*$ be the annihilation and creation operators 
in $L^2(\R^2)$, corresponding to the constant component
$B_0>0$ of the magnetic field $B$:
\begin{equation}
a_0=-2i e^{-B_0\abs{x}^2/4}\,\,\overline{\partial}\,e^{B_0\abs{x}^2/4},
\quad
a_0^*=-2i e^{B_0\abs{x}^2/4}\,\partial\, e^{-B_0\abs{x}^2/4}.
\label{e0}
\end{equation}
Then one can define 
\begin{equation}
a=e^{-\Psi_1}a_0e^{\Psi_1}
\text{ on }
\Dom(a)=\{e^{-\Psi_1}u\mid u\in\Dom(a_0)\}
\text{ and } a^*=e^{\Psi_1}a_0^*e^{-\Psi_1}.
\label{app1}
\end{equation}
In terms of these operators, we have 
\begin{equation}
h^+=a^*a \text{ and }
h^-=a a^*
\label{app2}
\end{equation}
and therefore $\s(h^+\setminus\{0\})=\s(h^-\setminus\{0\})$.
Finally, comparing the form ${\mathfrak h}^-$ with ${\mathfrak h}_0^-$,
one obtains
(see e.g. \cite{Besch} or \cite[Proposition 1.2]{Raikov3})
$$
{\mathfrak h}^-[u]\geq 
e^{2\essinf \Psi_1}{\mathfrak h}_0^-[e^{-\Psi_1}u]\geq 
2B_0e^{2\essinf \Psi_1}
\norm{e^{-\Psi_1}u}^2
\geq
2B_0 e^{-2\osc\Psi_1}\norm{u}^2,
\quad u\in\Dom({\mathfrak h}^-),
$$
where $\osc\Psi_1=\esssup\Psi_1-\essinf\Psi_1$.
It follows that $(0,m)$, $m=2B_0e^{-2\osc \Psi_1}>0$,
 is a gap in the spectrum of $h^+$ and of $h$.
See \cite{DubrovinN} for a different point of view on 
the issue of existence of the spectral gap 
and \cite{RozS} for further progress on this topic.

\textbf{3. Perturbations of the Pauli operator  and spectral asymptotics. }
Let $v\in L^p(\R^2)$, $p>1$, be a non-negative compactly supported
function, which has the physical meaning of the electric potential.
The Pauli operator which describes a particle in the external
magnetic field $B$ and electric field with the potential $\pm v$, is
\begin{equation}
h\pm v I= 
\begin{pmatrix}
h^+\pm v&0\\ 0 & h^-\pm v
\end{pmatrix}
\quad \text{in }L^2(\R^2)\oplus L^2(\R^2).
\label{a14}
\end{equation}
The main object of interest in this paper is the spectrum of $h^+\pm v$.
In order to define $h^+\pm v$ as a quadratic form sum, let us 
establish that $v$ is $h^+$-form compact. 
By \eqref{app1}, \eqref{app2} and boundedness of $\Psi_1$,
we see that $v$ is $h^+$-form compact if and only if
$ve^{-2\Psi_1}$ is $h_0^+$-form compact.
As $ve^{-2\Psi_1}\in L^p(\R^2)$, $p>1$, we obtain that 
$ve^{-2\Psi_1}$ is $h_0^+$-form compact
(see \cite{AHS}).

By the above established relative compactness,
the essential spectra of $h^++v$, $h^+-v$, and $h^+$ coincide;
moreover, due to the assumption $v\geq0$, the eigenvalues of $h^++v$
can accumulate to $0$ only from above, and the eigenvalues of 
$h^+-v$ can do so only from below.

Let $\lambda_1^-\leq\lambda_2^-\leq\cdots$ be the negative eigenvalues of
$h^+-v$, and $\lambda_1^+\geq\lambda_2^+\geq\cdots$ be the eigenvalues of 
$h^++v$ in the spectral gap $(0,m)$; 
here and in the rest of the paper, we assume
eigenvalues to be enumerated with multiplicities taken into account. 
The main aim of this paper is to describe the rate of convergence 
$\lambda_n^\pm\to0$ as $n\to\infty$. 
Roughly speaking, we prove the following asymptotics
(precise statements are given in Section~\ref{sec.b}):
\begin{equation}
\log(\pm n!\lambda_n^\pm)=n\log(B_0/2)+2n\log\Cp(\supp v)+o(n),
\quad n\to\infty,
\label{f6}
\end{equation}
where $\Cp$ is the \emph{logarithmic capacity} of a set.
The notion of logarithmic capacity is introduced in the framework
of potential theory; see e.g. \cite{Hille,Landkof}.
Recall that the logarithmic capacity of compact sets in $\R^2$
has the following properties:

(i) if $\Omega_1\subset\Omega_2$ then $\Cp\Omega_1\leq\Cp\Omega_2$;

(ii)
$\Cp\Omega$ coincides with the logarithmic capacity of the 
outer boundary of $\Omega$ ($=$the boundary of the unbounded
component of $\R^2\setminus\Omega$);

(iii)
the logarithmic capacity of a disc of radius $r$ is $r$;

(iv) if $\Omega_2=\{\alpha x\mid x\in\Omega_1\}$, $\a>0$,
then $\Cp \Omega_2=\a\Cp\Omega_1$.

We establish \eqref{f6} by means of the following 
simple chain of equivalent reformulations of the problem.
Firstly, a perturbation theory argument 
reduces the problem to the spectral asymptotics of an auxiliary 
compact self-adjoint operator $P_0vP_0$, 
where $P_0$ is the spectral projection of $h^+$, 
corresponding to the eigenvalue $0$.
Next, we observe that the eigenvalues of $P_0vP_0$
coincide with the singular numbers of a certain embedding 
operator (see \eqref{embed}).
Using the approach of \cite{Parfenov}, we relate the singular numbers
of this embedding operator to some sequence of orthogonal polynomials 
in the complex domain (see below).
Finally, application of the results of \cite{StahlTotik}
concerning the asymptotics of these orthogonal polynomials leads
to \eqref{f6}.

Using the same technique, we are also able to treat two similar
problems.  First, we consider the Pauli operators $h_0^+\pm v$ in the 
case of a constant magnetic field and describe the rate of accumulation
of the eigenvalues to the higher Landau levels. Secondly, we 
consider the three-dimensional Pauli Hamiltonian with a constant magnetic
field and a compactly supported electric potential and describe the 
rate of convergence of eigenvalues to 0. These results are presented in Section~2.

The rate of convergence of eigenvalues to zero for Pauli operators 
in dimensions two and three was investigated before 
in the case of constant magnetic field $B(x)=B_0>0$ for various classes 
of potentials $v$ with power or exponential decay at infinity;
see \cite{Solnyshkin,Sobolev,Tamura,Raikov1,Raikov2,Ivrii,RaikovWarzel}.
We refer the reader to the discussion in \cite{RaikovWarzel}.
The case of a constant magnetic field and 
compactly supported potentials $v$ was considered in 
\cite{RaikovWarzel} and \cite{MelgardRozenblum}.
The case of a two-dimensional operator with variable magnetic field and 
potentials $v$ with power or exponential  decay and also with compactly
supported potentials was treated in \cite{Raikov3}.
The results of \cite{RaikovWarzel,MelgardRozenblum,Raikov3}
for the case of compactly supported potentials read as
\begin{equation}
\log(\pm\lambda_{n}^{\pm})=-n\log n+O(n),
\quad n\to\infty.
\label{f7}
\end{equation}
As far as we are aware, a connection between the spectral asymptotics
of magnetic operators and logarithmic capacity or orthogonal
polynomials has not been made before.

Some physical intuition concerning these problem with constant
magnetic field can be gained from \cite{HornbergerS}.

\textbf{4. Orthogonal and Chebyshev polynomials.}
We identify $\R^2$ and $\C$ in a standard way:
$z=x_1+ix_2$ for $(x_1,x_2)\in\R^2$,
denote by $dm(z)$ the Lebesgue measure in $\C$ and 
consider $v$  as a function of $z$.
It appears that the sequence of polynomials in $z$,
orthogonal with respect to the measure $v(z)dm(z)$,
is related to the asymptotics of $\lambda_{n}^\pm$.
Here we present necessary facts from the theory of Chebyshev 
and orthogonal polynomials in complex domains; 
see e.g.  \cite{Hille} and \cite{StahlTotik} for the details.

For any $n=0,1,2,\dots$, let $\P_n$ be the set
of all monic polynomials in $z$ of degree $n$:
\begin{equation}
\P_n=\{z^n+a_{n-1}z^{n-1}+\dots+a_1z+a_0\mid a_0,\dots,a_{n-1}\in\C\}.
\label{a4}
\end{equation}
Let $\Omega\subset\C$ be a compact set.
For a fixed $n$, consider the problem of minimization of 
the norm $\norm{t}_{C(\Omega)}\equiv\sup_{z\in\Omega}\abs{t(z)}$
on the set $t\in\P_n$.
It is clear that the minimum is positive and attained at some polynomial 
$t_n\in\P_n$.
The polynomial $t_n$ is called the 
\emph{$n$'th Chebyshev polynomial} for the set $\Omega$.
(One can prove that such a polynomial is unique, but we will not 
need this fact).
It is well known that all zeros of $t_n$ lie in the closed convex hull of 
$\Omega$.
The \emph{$n$'th root asymptotics of $t_n$} is given by 
\begin{equation}
\lim_{n\to\infty}\norm{t_n}^{1/n}=\Cp\Omega.
\label{a2}
\end{equation}
Next, let $v\in L^1(\C,dm)$ be a non-negative compactly supported function.
Denote
\begin{equation}
M_n(v)=\inf_{p\in\P_n}\int_\C \abs{p(z)}^2 v(z)dm(z)=\int_\C \abs{p_n(z)}^2 v(z)dm(z),
\label{a5}
\end{equation}
where the sequence $\{p_n\}_{n=0}^\infty$, $p_n\in\P_n$, is obtained by applying the 
Gram--Schmidt orthogonalisation process in $L^2(\C, v(z)dm(z))$ to the 
sequence $1,z,z^2,\dots$. All zeros of $p_n$ lie in the closed convex hull 
of $\supp v$.

Regarding the $n$'th root asymptotics of $p_n$, the following 
facts are known (see \cite{StahlTotik}). Denote
\begin{equation}
\rho_+(v)=\limsup_{n\to\infty} M_n(v)^{1/n},
\quad
\rho_-(v)=\liminf_{n\to\infty} M_n(v)^{1/n}.
\label{a11}
\end{equation}
In general, it can happen that $\rho_-(v)<\rho_+(v)$
(see the proof of Theorem~1.1.9 in \cite{StahlTotik}).
One has the estimates 
\begin{gather}
\rho_+(v)\leq(\Cp\supp v)^2,
\label{a12}
\\
\rho_-(v)\geq
(\Cp\Omega_-(v))^2, \quad
\Omega_-(v)=
\{z\in\C\mid\limsup_{r\to+0}
\frac{\log\int_{\abs{z-\zeta}\leq r}v(\zeta)dm(\zeta)}{\log r}<\infty\}.
\label{a13}
\end{gather}
The inequality \eqref{a12} is a part of Corollary~1.1.7 of \cite{StahlTotik}.
The inequality \eqref{a13}, although not stated explicitly in \cite{StahlTotik},
follows directly from the proof of Theorem~4.2.1 therein.

\begin{remark}\label{rmk1}
Let $\Omega\subset\C$ be a compact set with a Lipschitz 
boundary, and let $v\in L^1(\C,dm)$ be such that $v(z)\geq c>0$ for all 
$z\in\Omega$ and $v(z)=0$ for all $z\in\C\setminus\Omega$.
Then we easily find that 
$\Omega_-(v)=\Omega=\supp v$
and therefore $\rho_+(v)=\rho_-(v)=(\Cp\Omega)^2$.
\end{remark}

%%%%%%%%%%%%%%%%%%%%%%%%%%%%%%%%%%%%%%%%
\section{Main results}\label{sec.b}
%%%%%%%%%%%%%%%%%%%%%%%%%%%%%%%%%%%%%%%%
\textbf{1. Two-dimensional Pauli operators with variable magnetic field.}
Let $h^+$, as in the Introduction, be the Pauli operator defined 
via \eqref{f4} with $\Psi$ subject to \eqref{psi}.
Let $v$ and $\lambda_n^\pm$ be as in the Introduction.
\begin{theorem}\label{th1}
Let $0\leq v\in L^p(\R^2)$, $p>1$, be a compactly supported potential
and let $M_n(v)$ be as defined in \eqref{a5}.
Then there exists $k\in\N$ such that 
\begin{align}
(B_0/2)M_{n+k}(v)^{1/n}(1+o(1))
\leq (n!\lambda^+_{n})^{1/n}\leq
(B_0/2)M_{n-1}(v)^{1/n}(1+o(1)),
\label{a6a}
\\
(B_0/2)M_{n-1}(v)^{1/n}(1+o(1))
\leq (-n!\lambda^-_{n})^{1/n}\leq
(B_0/2)M_{n-k}(v)^{1/n}(1+o(1)),
\label{a6b}
\end{align}
as $n\to\infty$.
In particular, 
$$
\limsup_{n\to\infty}(\pm n!\lambda_{n}^\pm)^{1/n}=B_0\rho_+(v)/2,
\quad
\liminf_{n\to\infty}(\pm n!\lambda_{n}^\pm)^{1/n}=B_0\rho_-(v)/2,
$$
where $\rho_\pm(v)$ are defined by \eqref{a11}.
If $v$ is of the class described in Remark~\ref{rmk1},
then the asymptotics \eqref{f6} holds true.
\end{theorem}
\begin{remark}
Let $\mu$ be a compactly supported finite measure in $\R^2$ 
such that the quadratic form $\int_\C\abs{u(x)}^2d\mu(x)$
is compact with respect to the quadratic form 
${\mathfrak h}^+$. Then one can define the 
self-adjoint operators  corresponding to the 
quadratic forms 
$$
{\mathfrak h}^+[u]\pm\int_{\R^2}\abs{u(x)}^2d\mu(x).
$$
All our considerations remain valid for such operators.
For example, the case of a measure $\mu$, supported 
by a curve, can be interesting.
\end{remark}

\textbf{2. Two-dimensional Pauli operators with constant magnetic field.}
Let $B(x)=B_0>0$; consider the corresponding 
operator  $h_0^+$.
As it is well known, the spectrum of $h_0^+$ consists of the eigenvalues 
$\{2qB_0\}_{q=0}^\infty$ of infinite multiplicities; 
these eigenvalues are known as Landau levels.
Consider the problem of accumulation of eigenvalues
of $h_0^+\pm v$ to a fixed higher Landau level $2qB_0$, $q\geq1$.
Let $\lambda_{q,1}^-\leq\lambda_{q,2}^-\leq\cdots$ be the eigenvalues of $h_0^+-v$
in the interval $(2(q-1)B_0,2qB_0)$, and let $\lambda_{q,1}^+\geq\lambda_{q,2}^+\geq\cdots$
be the eigenvalues of $h_0^++v$ in $(2qB_0,2(q+1)B_0)$.
\begin{theorem}\label{th2}
Let $\Omega\subset\R^2$ be a compact set with Lipschitz boundary
and let $v\in L^p(\R^2)$, $p>1$, be such that $v(x)\geq c>0$ for 
$x\in\Omega$ and $v(x)=0$ for $x\in\R^2\setminus\Omega$.
Then for the corresponding eigenvalues $\lambda_{q,n}^\pm$ we have:
$$
\lim_{n\to\infty}(\pm n!(\lambda_{q,n}^\pm-2qB_0))^{1/n}=\frac{B_0}{2}(\Cp\Omega)^2.
$$
\end{theorem}
The  rate of convergence of $\lambda_{q,n}^\pm \to2qB_0$ as $n\to\infty$, $q\geq1$, 
was studied before in \cite{RaikovWarzel,MelgardRozenblum}, 
where the asymptotics
$$
\log(\pm(\lambda_{q,n}^\pm-2qB_0))=-n\log n+O(n),
\quad n\to\infty
$$
was obtained. Note that if the potential $v$ depends only on $\abs{x}$,
then the result of Theorem~\ref{th2} can be obtained by a direct calculation
using separation of variables, see e.g. \cite[Proposition 3.2]{RaikovWarzel}.

\textbf{3. Three-dimensional Pauli operator with a constant magnetic field.} 
Let 
$$
H=(-i\nabla-\mathbf{A}(\mathbf{x}))^2-B_0
\text{ in }L^2(\R^3,d\mathbf{x}),
\quad \mathbf{x}=(x_1,x_2,x_3),
$$
where $\mathbf{A}(\mathbf{x})=(-\frac12 B_0x_2,\frac12 B_0x_1,0)$.
It is well known that the spectrum of $H$ is absolutely continuous 
and coincides with the interval $[0,\infty)$. 
The background information concerning the spectral theory of $H$ and its perturbations
can be found in \cite{AHS}.

Let $V\in L^{3/2}(\R^3)$ be a non-negative compactly supported potential.
The operator of multiplication by $V$ in $L^2(\R^3)$ is $H$-form compact
(cf. \cite{AHS}).
Thus, one can define the self-adjoint opertor $H-V$ via the corresponding 
quadratic form; the essential spectrum of $H-V$ is also $[0,\infty)$.
Let $\Lambda_1\leq\Lambda_2\leq\cdots$ be the negative eigenvalues of $H-V$;
we have $\Lambda_n\to0$ as $n\to\infty$.
Below we describe the asymptotic behaviour of $\Lambda_n$ as $n\to\infty$
in terms of the auxiliary weight function 
\begin{equation}
w(x_1,x_2)=\int_{-\infty}^\infty V(x_1,x_2,x_3)dx_3,
\quad \text{a.e. } (x_1,x_2)\in\R^2.
\label{a0}
\end{equation}
As above, we consider $w$ as a function of $z=x_1+ix_2$.
\begin{theorem}\label{th3}
Let $0\leq V\in L^{3/2}(\R^3)$ be a compactly supported potential
and $w$ be defined by \eqref{a0}.
Then there exists $k\in\N$ such that 
\begin{equation}
(B_0/2)^2M_{n+k}(w)^{2/n}(1+o(1))
\leq (-(n!)^2\Lambda_n)^{1/n}\leq
(B_0/2)^2M_{n-k}(w)^{2/n}(1+o(1)),
\label{a1}
\end{equation}
as $n\to\infty$.
In particular, 
$$
\limsup_{n\to\infty}(-(n!)^2\Lambda_n)^{1/n}=(B_0\rho_+(w)/2)^2,
\quad
\liminf_{n\to\infty}(-(n!)^2\Lambda_n)^{1/n}=(B_0\rho_-(w)/2)^2,
$$
where $\rho_\pm(w)$ are defined by \eqref{a11}.
\end{theorem}

The rate of accumulation $\Lambda_n\to0$ for potentials $V$ with power
or exponential decay was considered before in 
\cite{Solnyshkin,Sobolev,Tamura,Raikov1,Raikov2,Ivrii,RaikovWarzel}.
For compactly supported potentials, this problem was considered 
in \cite{RaikovWarzel,MelgardRozenblum}, where the asymptotics
$$
\log \Lambda_n=-2n\log n+O(n),
\quad n\to\infty
$$
was obtained.

\begin{remark}
Theorem~\ref{th3} remains valid under the following
assumptions on $V$:
(i) $V\geq0$, $V$ is $H$-form compact; 
(ii) $\int_{\R^3} V(\mathbf{x})(1+\abs{x_3}^2)d\mathbf{x}<\infty$;
(iii) the  function $w$, defined by \eqref{a0}, is
compactly supported.
\end{remark}

%%%%%%%%%%%%%%%%%%%%%%%%%%%%%%%%%%%%
\section{Proof of Theorems~\protect\ref{th1}, \protect\ref{th2}
and \protect\ref{th3}}
\label{sec.c}
%%%%%%%%%%%%%%%%%%%%%%%%%%%%%%%%%%%%
\emph{Proof of Theorem~\ref{th1}:}
Let $\H_0\subset L^2(\R^2)$ be the kernel of $h^+$, 
and let $P_0$ be the corresponding eigenprojection,
$\Ran P_0=\H_0$.
Consider the compact self-adjoint operator $P_0vP_0$.
The key ingredient in the proof is the following 
\begin{lemma}\label{lma1}
Let $v\in L^1(\R^2)$ be a non-negative compactly supported function
and let $s_1\geq s_2\geq\cdots>0$
be the eigenvalues of $P_0vP_0$. Then
\begin{equation}
(n! s_{n+1})^{1/n}=
(B_0/2)M_n(v)^{1/n}(1+o(1)),
\quad n\to\infty,
\label{c3}
\end{equation}
where $M_n(v)$ are defined by \eqref{a5}.
\end{lemma}
The proof is given in Section~\ref{sec.d}.
Now it remains to employ a perturbation theory argument
(see \cite[Proposition 3.1]{Raikov3}
or \cite[Proposition 4.1]{RaikovWarzel}) 
based on the Birman-Schwinger principle and on 
Weyl inequalities for eigenvalues of a sum of 
compact operators.
This argument shows that there exists $k\in\N$ 
such that for all sufficiently large $n\in\N$
one has 
$$
s_{n}\leq -\lambda_n^-\leq 2 s_{n-k}.
\quad
\frac12 s_{n+k}\leq \lambda_n^+\leq s_{n}.
$$
Combining these inequalities with 
Lemma~\ref{lma1},
we obtain the required result.
\qed

\emph{Proof of Theorem~\ref{th2}:}
For any $q\geq0$, denote $\H_q=\Ker(h^+-2qB_0)$ and let 
$P_q$ be the eigenprojection of $h_0^+$ corresponding to the 
eigenvalue $2qB_0$. 
Consider the compact self-adjoint operator $P_qvP_q$, and let 
$s_1^{(q)}\geq s_2^{(q)}\geq\cdots$
be the eigenvalues of this operator. 
As in the proof of Theorem~\ref{th1},
using a perturbation theory argument 
based on the Birman-Schwinger principle
and Weyl inequalities 
(see \cite[Proposition~4.1]{RaikovWarzel}),
one shows that there exists $k\in\N$ such that for all 
sufficiently large $n\in\N$,
\begin{equation}
\frac12 s^{(q)}_{n+k}\leq\pm(\lambda_{q,n}^\pm-2qB)\leq 2 s^{(q)}_{n-k}.
\label{e5}
\end{equation}

Now the proof of Theorem~\ref{th2} reduces to 
\begin{lemma}\label{lma2}
Let $\Omega\subset\R^2$ be a compact set with Lipschitz boundary
and let $v\in L^1(\R^2)$ be such that $v(x)\geq c>0$
for $x\in\Omega$ and $v(x)=0$ for $x\in\R^2\setminus\Omega$.
Fix $q\in\N$ and let 
$s^{(q)}_1\geq s^{(q)}_2\geq\cdots$
be the eigenvalues of $P_qvP_q$.
Then one has
\begin{equation}
\lim_{n\to\infty}(n! s_n^{(q)})^{1/n}=(B_0/2)(\Cp \Omega)^2.
\label{c6}
\end{equation}
\end{lemma}

The proof of Lemma~\ref{lma2} 
is given in Section~\ref{sec.e}. From Lemma~\ref{lma2}
and the estimate \eqref{e5}, we immediately obtain the required result.
\qed

\emph{Proof of Theorem~\ref{th3}:}
The proof repeats almost word for word the construction of \cite{Sobolev}.
According to the Birman-Schwinger principle, for $E>0$ we have:
\begin{equation}
\sharp\{n\mid \Lambda_n<-E\}=n_+(1;\sqrt{V}(H_0+E)^{-1}\sqrt{V}).
\label{e1}
\end{equation}
The operator $\sqrt{V}(H_0+E)^{-1}\sqrt{V}$ can be represented as
$$
\sqrt{V}(H_0+E)^{-1}\sqrt{V}=\frac{1}{2\sqrt{E}}K_1+K_2+K_3.
$$
Here $K_1$, $K_2$ are the operators in $L^2(\R^3)$ with the 
integral kernels
\begin{align*}
K_1(\mathbf{x},\mathbf{y})
&=
\sqrt{V(\mathbf{x})}P_0(x_\perp,y_\perp)\sqrt{V(\mathbf{y})},
\\
K_2(\mathbf{x},\mathbf{y})
&=
\sqrt{V(\mathbf{x})}P_0(x_\perp,y_\perp)
\frac{e^{-\sqrt{E}\abs{x_3-y_3}}-1}{2\sqrt{E}}
\sqrt{V(\mathbf{y})},
\end{align*}
where the notation $x_\perp=(x_1,x_2)$, $y_\perp=(y_1,y_2)$ is used,
and $P_0(x_\perp,y_\perp)$ is the integral kernel of the 
operator $P_0$ in $L^2(\R^2)$.
Finally, $K_3$ is the operator 
$$
K_3=\sqrt{V}Q_0(H_0+E)^{-1}\sqrt{V},
$$
where $Q_0=(I-P_0)\otimes I$
in the decomposition
$L^2(\R^3,dx_1 dx_2 dx_3)=L^2(\R^2,dx_1 dx_2)\otimes L^2(\R,dx_3)$.

The operators $K_2$ and $K_3$ have limits (in the operator norm)
as $E\to+0$; these limits are compact self-adjoint operators.
Thus, by the Weyl's inequalities for eigenvalues
(see e.g. \cite{BirmanS}), we have 
for $E\to+0$:
\begin{gather}
n_+(1;\sqrt{V}(H_0+E)^{-1}\sqrt{V})
\leq 
n_+(\frac12;\frac{1}{2\sqrt{E}}K_1)+n_+(\frac12;K_2+K_3)
\leq
n_+(\sqrt{E};K_1)+O(1),
\label{e2}
\\
n_+(1;\sqrt{V}(H_0+E)^{-1}\sqrt{V})
\geq 
n_+(\frac32;\frac{1}{2\sqrt{E}}K_1)-n_+(\frac12;-K_2-K_3)
\geq
n_+(3\sqrt{E}; K_1)-O(1).
\label{e3}
\end{gather}
Finally, again as in \cite{Sobolev}, let us prove that the non-zero
eigenvalues of $K_1$ coincide with those of $P_0wP_0$, where
$w$ is defined by \eqref{a0}.
It suffices to prove this statement for continuous $V$ with compact support;
the general case $V\in L^{3/2}$ then follows by approximation argument.
Let $N_1: L^2(\R^3,dx_1\,dx_2\,dx_3)\to L^2(\R^2,dx_1\,dx_2)$ and $N_2: L^2(\R^2,dx_1\,dx_2)\to  L^2(\R^3,dx_1\,dx_2\,dx_3)$ be the following operators:
\begin{gather*}
(N_1u)(x_1,x_2)=\int_{-\infty}^\infty V^{1/2}(x_1,x_2,x_3)u(x_1,x_2,x_3)dx_3,
\\
(N_2u)(x_1,x_2,x_3)=V^{1/2}(x_1,x_2,x_3)u(x_1,x_2).
\end{gather*}
Then $K_1=N_2P_0N_1=(N_2P_0)(P_0N_1)$ and 
$P_0wP_0=(P_0N_1)(N_2P_0)$.
It follows that the non-zero eigenvalues of $K_1$ coincide
with $\{s_n\}$, the non-zero eigenvalues of $P_0wP_0$,
and so $n_+(\sqrt{E};K_1)=\sharp\{n\mid (s_n)^2>E\}$.
{}From here and \eqref{e1}, \eqref{e2}, \eqref{e3}
it follows that for some $k\in\N$ and all sufficiently large $n\in\N$, 
one has
\begin{equation}
\frac19(s_{n+k})^2\leq \Lambda_n\leq (s_{n-k})^2.
\label{e4}
\end{equation}
Combining this with Lemma~\ref{lma1}, 
we get the statement of Theorem~\ref{th3}.
\qed

%%%%%%%%%%%%%%%%%%%%%%%%%%%%%%%%%%%%
\section{Proof of Lemma~\protect\ref{lma1}}\label{sec.d}
%%%%%%%%%%%%%%%%%%%%%%%%%%%%%%%%%%%%
First let us consider the case of a constant magnetic field
$B(x)=B_0>0$.
Let $F^2$ be the Hilbert space of all entire functions 
$f$ such that
$$
\norm{f}^2_{F^2}
=\int_\C\abs{f(z)}^2 e^{-B_0\abs{z}^2/2}dm(z)<\infty.
$$
In the case $B_0=2$, the space $F^2$ 
is usually called Fock space or Segal-Bargmann space.
By \eqref{a16}, we have an isometry between 
$\H_0=\Ker h^+\subset L^2(\C,dm)$ and $F^2$, given by 
$u(z)=e^{-B_0\abs{z}^2/4} f(z)$, $u\in \H_0$, $f\in F^2$.
Thus, the quadratic form of the operator $P_0 vP_0\mid_{\H_0}$
is unitarily equivalent to the quadratic form 
$$
\int_\C \abs{f(z)}^2 v(z)e^{-B_0\abs{z}^2/2}dm(z),
\quad f\in F^2.
$$
It follows that the non-zero eigenvalues $s_n$ of $P_0vP_0$
coincide with the singular values $\mu_n$ of the embedding
operator 
\begin{equation}
F^2\subset L^2(\C,v(z)e^{-B_0\abs{z}^2/2}dm(z)).
\label{embed}
\end{equation}
The case of a variable magnetic field can be also reduced to the 
embedding \eqref{embed}.
Indeed, using the boundedness of $\Psi_1$, one obtains
(see \cite[Proposition 3.2]{Raikov3}):
$$
\mu_n
e^{-2\osc\Psi_1}
\leq
s_n
\leq
\mu_n e^{2\osc\Psi_1},
\quad n\in\N.
$$
Thus, it remains to prove  the asymptotic formula
\begin{equation}
(n!\mu_{n+1})^{1/n}=(B_0/2)M_n(v)^{1/n}(1+o(1)),
\quad n\to\infty
\label{c13}
\end{equation}
for the singular values $\mu_n$ of the embedding \eqref{embed}.
We shall assume $B_0=2$; the general case can be reduced 
to this one by a linear change of coordinates.

Asymptotics of the $n$-widths of the embedding 
$F^2\subset C(\Omega)$, where $\Omega$ is a compact set in $\C$, was
studied in \cite{Parfenov}. Below we repeat the arguments of 
\cite{Parfenov} (with trivial modifications) to obtain the required asymptotics.

By the minimax principle, we have 
 the following variational characterisation of $\mu_n$:
\begin{align}
\mu_{n+1}
&=
\inf_{L_n^+\subset F^2}\sup_{f\in L_n^+\setminus\{0\}}
\frac{\int_{\C}\abs{f(z)}^2 v(z)e^{-\abs{z}^2}dm(z)}{\norm{f}^2_{F^2}},
\quad\codim L_n^+=n,
\label{b6}
\\
\mu_{n+1}
&=
\sup_{L_n^-\subset F^2}\inf_{f\in L_n^-\setminus\{0\}}
\frac{\int_{\C}\abs{f(z)}^2 v(z)e^{-\abs{z}^2}dm(z)}{\norm{f}^2_{F^2}},
\quad\dim L_n^-=n+1.
\label{b5}
\end{align}

\textbf{1. Upper bound on $\mu_{n+1}$.}
For the subspaces $L_n^+$ from \eqref{b6}, we will take 
$$
L_n^+=\{f\in F^2\mid f(z)=p_n(z)g(z), \text{ $g$ is entire function}\},
$$
where $p_n$ is the sequence of monic  polynomials
orthogonal with respect to the measure  
$v(z)dm(z)$.

In order to estimate the ratio 
in \eqref{b6} from above, let us prove the following auxiliary 
statement.
Denote $R_0=\max_{z\in\supp v}\abs{z}$.
We claim that for any $\e\in(0,\frac13)$, there exists 
$N\in\N$ such that for all $n\geq N$ and any $f=p_n g\in L_n^+$,
we have 
\begin{equation}
\sup_{\abs{z}\leq R_0}\abs{g(z)}^2\leq (1-\e)^{-2n}
\frac{1}{n!}\norm{p_n g}^2_{F^2}.
\label{c9}
\end{equation}
Indeed, we have 
$$
g(z)
=
\frac1{2\pi i}\int_{\abs{\zeta}=r}\frac{f(\zeta)}{p_n(\zeta)(\zeta-z)}d\zeta,
\quad r>R_0,
$$
and therefore
$$
\sup_{\abs{z}\leq R_0}\abs{g(z)}^2
\leq
\frac{r}{2\pi}\sup_{\abs{z}\leq R_0}
\int_{\abs{\zeta}=r}
\frac{\abs{f(\zeta)}^2}{\abs{p_n(\zeta)}^2\abs{\zeta-z}^2}d\abs{\zeta}
$$
for any $r>R_0$.
Denote $R=R_0/\e$.
Since all zeros of $p_n$ lie in the closed convex hull
of $\supp v$, we obtain:
$$
\abs{p_n(\zeta)}\abs{\zeta-z}\geq ((1-\e)r)^{n+1}, 
\quad \abs{z}\leq R_0, \quad \abs{\zeta}=r\geq R.
$$
Thus, we get 
$$
\sup_{\abs{z}\leq R_0}\abs{g(z)}^2
\leq
\frac{r^{-2n-1}}{2\pi(1-\e)^{2n+2}}
\int_{\abs{\zeta}=r}\abs{f(\zeta)}^2d\abs{\zeta},
\quad r\geq R.
$$
Integrating the last inequality over $r$ from $R$ to $\infty$ 
with the weight $e^{-r^2}r^{2n+1}$, 
and using the fact that 
$$
\int_R^\infty e^{-r^2}r^{2n+1}dr
=
\frac12 n!
-
\int_0^R e^{-r^2}r^{2n+1}dr
\geq
\frac1{2\pi}(1-\e)^{-2} n!
$$
for all sufficiently large $n$, we obtain \eqref{c9}.

From \eqref{c9} we obtain for any $f=p_n g\in L_n^+$:
$$
\int_\C\abs{f(z)}^2v(z)e^{-\abs{z}^2}dm(z)
\leq
M_n(v)\norm{g}^2_{C(\supp v)}
\leq
M_n(v) (1-\e)^{-2n}\frac1{n!}\norm{f}_{F^2}^2.
$$
Together with \eqref{b6}, the last estimate yields
\begin{equation}
(n! \mu_{n+1})^{1/n}\leq (1-\e)^{-2} M_n(v)^{1/n}
\label{c10}
\end{equation}
for all sufficiently large $n$.

\textbf{2. Lower bound for $\mu_{n+1}$.}
Let us use formula \eqref{b5} and take $L_n^-$ to be the set 
of all polynomials in $z$ of degree $\leq n$.
As in the proof of the upper bound, we denote 
$R_0=\max_{z\in\supp v}\abs{z}$, fix $\e>0$ and set 
$R=R_0/\e$.
We shall use the following equivalent norm in $F^2$:
\begin{equation}
\Norm{f}_{F^2}^2=\int_{\abs{z}\geq R}\abs{f(z)}^2 e^{-\abs{z}^2} dm(z),
\quad 
\Norm{f}_{F^2}\leq \norm{f}_{F^2}\leq C(R)\Norm{f}_{F^2}.
\label{d1}
\end{equation}
Let $q_n\in L_n^-\setminus\{0\}$ be the polynomial which 
minimizes the ratio 
\begin{equation}
\frac{\int_\C\abs{q_n(z)}^2v(z)e^{-\abs{z}^2}dm(z)}{\Norm{q_n}^2_{F^2}}
\label{d2}
\end{equation}
among all polynomials in $L_n^-\setminus\{0\}$.
The following standard argument shows that all zeros of $q_n$ 
are confined to the disk $\{z\mid \abs{z}\leq R_0\}$.
Suppose that one of the zeros $z_k$ is outside the disk; 
then replace $q_n(z)$ by 
$q_n(z)\abs{z_k}(z-R_0^2/\overline{z}_k)/(R_0(z-z_k))$.
One has 
$$
\frac{\abs{z_k}\abs{z-R_0^2/\overline{z}_k}}{R_0\abs{z-z_k}}\leq1
\text{ for $\abs{z}\leq R_0$ and }
\frac{\abs{z_k}\abs{z-R_0^2/\overline{z}_k}}{R_0\abs{z-z_k}}\geq1
\text{ for $\abs{z}\geq R_0$,}
$$
so this change 
decreases the ratio \eqref{d2} --- contradiction. 
Next, without the loss of generality, we may assume that $q_n$ is 
monic.
Denote $m=\deg q_n\leq n$; we get the estimate 
$$
\Norm{q_n}_{F^2}^2
=
\int_{\abs{z}\geq R}\abs{q_n(z)}^2 e^{-\abs{z}^2} dm(z)
\leq
\int_{\abs{z}\geq R}\abs{z}^{2m}(1+\e)^{2m}e^{-\abs{z}^2} dm(z)
\leq
(1+\e)^{2m}\pi m!.
$$
On the other hand, for the numerator of \eqref{d2}, we have 
$$
\int_{\C}\abs{q_n(z)}^2 v(z)e^{-\abs{z}^2}dm(z)
\geq
e^{-R_0^2}\int_{\C} \abs{q_n(z)}^2v(z)dm(z)
\geq
e^{-R_0^2}M_m(v).
$$
Combining the above estimates, we obtain:
\begin{equation}
\mu_{n+1}
\geq
\inf_{f\in L_n^-\setminus\{0\} }
\frac{\int_{\C}\abs{f(z)}^2 v(z)e^{-\abs{z}^2}dm(z)}{C(R)\Norm{f}^{2}_{F^2}}
\geq
\min_{0\leq m\leq n}
\frac{M_m(v)}{C_1(R)(1+\e)^{2m} m!}.
\label{d3}
\end{equation}
As $zp_m(z)\in\P_{m+1}$, from the definition \eqref{a5} of $M_m(v)$ 
we get a trivial estimate $M_{m+1}(v)\leq R_0^2 M_m(v)$.
This estimate shows that for a sufficiently large $n$, the minimum 
in \eqref{d3} is attained at $m=n$.
Therefore,
$$
(n! \mu_{n+1})^{1/n}
\geq
\left(\frac{1}{C_1(R)}\right)^{1/n}
\frac{M_n(v)^{1/n}}{(1+\e)^2}
\geq
(1+\e)^{-3}M_n(v)^{1/n}
$$
for all sufficiently large $n$.
The latter estimate together with \eqref{c10} completes the proof of
the Lemma.
\qed

%%%%%%%%%%%%%%%%%%%%%%%%%%%%%%%

\section{Proof of Lemma~\protect\ref{lma2}}\label{sec.e}
%%%%%%%%%%%%%%%%%%%%%%%%%%%%%%%
First recall some well known facts concerning the spectral decomposition 
of the operator $h_0^+$.
The operator  $h_0^+$ can be represented in terms of 
 the annihilation and creation operators \eqref{e0}
as $h_0^+=a_0^*a_0$.
The operators $a_0$, $a_0^*$ obey the commutation relation $[a_0,a_0^*]=2B_0$,
wherefrom we get the identity
$a_0^q(a_0^*)^q u=(2B_0)^q q! u$ for all 
$u\in\H_0$ and $q\in\N$. 
It follows that 
\begin{equation}
(2B_0)^{-q/2}(q!)^{-1/2}(a_0^*)^q:\H_0\to\H_q
\text{ is an isometry onto $\H_q$.}
\label{c2}
\end{equation}
Recalling the explicit isomorphism between $\H_0$ and the space 
$F^2$ (see the previous section), we see that the change 
$u=(2B_0)^{-q/2}(q!)^{-1/2}(a_0^*)^q(e^{-B_0\abs{z}^2/2}f(z))$
gives a unitary equivalence between the operator $P_qvP_q$
and the operator in $F^2$ defined by the quadratic form
\begin{equation}
(2B_0)^{-q}(q!)^{-1}\int_\C\abs{(a_0^*)^q(e^{-B_0\abs{z}^2/2}f(z))}^2v(z)dm(z).
\label{d4}
\end{equation}
For simplicity we will consider the case $q=1$; the general case can be treated
in a similar manner.
Also, we will take $B_0=2$; the general case can be reduced to this one 
by a linear change of variables.
With these simplifications, the form \eqref{d4} becomes
\begin{equation}
\int_\C\abs{f'(z)-\overline{z}f(z)}^2 e^{-\abs{z}^2}v(z)dm(z),
\quad f\in F^2.
\label{c7}
\end{equation}
Let us prove the asymptotics \eqref{c6} for the eigenvalues
$\{s_n^{(1)}\}_{n=1}^\infty$ corresponding to the form \eqref{c7}.

\textbf{1. Upper bound for $s_n^{(1)}$.}
Let $\Omega_\delta=\{z\mid \dist(z,\Omega)\leq \delta\}$.
By the Cauchy integral formula, we have
$$
\sup_{\Omega}\abs{f'}\leq\frac1\delta\sup_{\Omega_\delta}\abs{f},
\quad f\in F^2.
$$
Thus, we have the following bound for the form \eqref{c7}:
$$
\int_\C\abs{f'(z)-\overline{z}f(z)}^2 e^{-\abs{z}^2}v(z)dm(z)
\leq C\norm{f}^2_{C(\Omega_\d)},
$$
where $C$ depends on $\Omega$, $\delta$, $v$.
Let us define 
$$
L_n^+=\{f\in F^2\mid f(z)=t_n(z)g(z), \text{$g$ entire}\},
$$ 
where $t_n$ is the $n$'th Chebyshev polynomial for the set
$\Omega_\delta$.
Note that the proof of \eqref{c9} uses only the fact that 
all zeros of $p_n$ lie in the closed convex hull
of $\supp v$. 
Therefore, the same estimate remains valid with the change 
$p_n\mapsto t_n$.
Thus, we get
\begin{multline*}
\int_{\C}\abs{f'(z)-\overline{z}f(z)}^2 
e^{-\abs{z}^2}v(z)dm(z)
\leq
C\norm{f}^2_{C(\Omega_\d)}
\leq
C\norm{t_n}^2_{C(\Omega_\delta)}
\norm{g}^2_{C(\Omega_\delta)}
\\
\leq
C(1-\e)^{-2n}\frac1{n!}
\norm{t_n}^2_{C(\Omega_\delta)}\norm{f}^2_{F^2},
\quad f\in L_n^+,
\end{multline*}
for all sufficiently large $n$.
This yields
$$
\limsup_{n\to\infty}(n!s_n^{(1)})^{1/n}
\leq
(1-\e)^{-2}
\lim_{n\to\infty}\norm{t_n}_{C(\Omega_\d)}^{2/n}
=
(1-\e)^{-2}
(\Cp\Omega_\d)^2.
$$
It remains to note that $\e$ and $\delta$
can be chosen arbitrary small and that
$\lim_{\d\to+0}\Cp\Omega_\d=\Cp\Omega$
for any compact set $\Omega$ (see e.g. \cite{Hille}).

\textbf{2. Lower bound for $s_{n}^{(1)}$.}
Due to the compactness of the embedding of the Sobolev space 
$W_2^1(\Omega)\subset L^2(\Omega)$,
for any $\gamma>0$
there exists a subspace $N\subset W^1_2(\Omega)$
of a finite codimension such that 
\begin{equation}
\norm{u}_{L^2(\Omega)}\leq 
\gamma\norm{\nabla u}_{L^2(\Omega)},
\quad \forall u\in N.
\label{c11}
\end{equation}
Let us choose $\gamma=1/(4R_0)$, 
$R_0=\max_{\Omega}\abs{z}$,
consider the corresponding subspace $N$ and 
denote $\codim N=l<\infty$.
Next, let $L_n^-$, as above, be the set of all polynomials in $z$ 
of degree $\leq n$.
Consider the subspace $\tilde L_n^-=L_n^-\cap N$;
clearly, $\dim \tilde L_n^-\geq n+1-l$.
By \eqref{c11}, for any $f\in \tilde L_n^-$ we have 
$\norm{f}_{L^2(\Omega)}\leq \frac1{2R_0}\norm{f'}_{L^2(\Omega)}$.
It follows that 
for any $f\in \tilde L_n^-$:
$$
\norm{f'}_{L^2(\Omega)}
\leq 
\norm{f'-\overline{z}f}_{L^2(\Omega)}+\norm{\overline{z}f}_{L^2(\Omega)}
\leq
\norm{f'-\overline{z}f}_{L^2(\Omega)}+\frac12\norm{f'}_{L^2(\Omega)},
$$
and so 
$$
\norm{f'-\overline{z}f}_{L^2(\Omega)}\geq
\frac12\norm{f'}_{L^2(\Omega)}\geq
R_0\norm{f}_{L^2(\Omega)}.
$$
Thus, for the quadratic form \eqref{c7} we have 
$$
\int_\C\abs{f'(z)-\overline{z}f(z)}^2 e^{-\abs{z}^2}v(z)dm(z)
\geq
C\int_{\Omega}\abs{f(z)}^2dm(z),
\quad f\in \tilde L_n^-.
$$
According to the second inequality in \eqref{d3} 
with $v=\chi_\Omega$
(we denote by $\chi_\Omega$ the characteristic function of $\Omega$)
we have therefore
$$
\frac{\int_\C\abs{f'(z)-\overline{z}f(z)}^2 
e^{-\abs{z}^2}v(z)dm(z)}{\norm{f}^2_{F^2}}
\geq
C\frac{\int_\Omega\abs{f(z)}^2dm(z)}{\norm{f}^2_{F^2}}
\geq
C'
\frac{M_{n}(\chi_\Omega)}{(1+\e)^{2n}n!},
\quad 
f\in\tilde L_n^-.
$$
It follows that 
$$
s^{(1)}_{n+1-l}\geq C'\frac{M_{n}(\chi_\Omega)}{(1+\e)^{2n}n!},
$$
for all sufficiently large $n$.
As stated in Remark~\ref{rmk1},
$\lim_{n\to\infty} M_n(\chi_\Omega)^{1/n}=(\Cp \Omega)^2$.
Thus, 
$$
\liminf_{n\to\infty}(n! s_n^{(1)})^{1/n}\geq\frac{(\Cp\Omega)^2}{(1+\e)^2}
$$
for any $\e>0$.
\qed

\section*{Acknowledgements}
We are indebted to  M.~Sh.~Birman, Yu.~Netrusov, G.~Raikov, 
G.~Rozenblum, A.~Sobolev and H.~Stahl   for useful discussions.
The work was supported by the Royal Society grant 2004/R1-FS.
The authors are grateful to the Mathematisches Forschungsinstitut Oberwolfach 
for hospitality and financial support. The first named author is also grateful to 
Loughborough University for hospitality.


\begin{thebibliography}{10}

\bibitem{AharonovC}
Y.~Aharonov, A.~Casher, 
{\it Ground state of a spin-${\frac12}$ charged particle in a two-dimensional magnetic field.} 
 Phys. Rev. A (3) {\bf 19} (1979), no.~6, 2461--2462.


\bibitem{AHS}
J.~Avron, I.~Herbst\ and\ B.~Simon, 
{\it Schr\"odinger operators with magnetic fields. I. General interactions},
Duke Math. J. {\bf 45} (1978), no.~4, 847--883.


\bibitem{Besch}
A.~Besch,
{\it Eigenvalues in spectral gaps of the two-dimensional Pauli operator.}
J. Math. Phys. {\bf 41}, 7918--7931.

\bibitem{BirmanS}
M.~Sh.~Birman and M.~Z.~Solomyak,
{\it Spectral theory of self-adjoint operators in Hilbert space,}
Dordrecht, D.Reidel P.C., 1987

\bibitem{DubrovinN}
B.~A.~Dubrovin, S.~P.~Novikov, 
{\it Fundamental states in a periodic field. 
Magnetic Bloch functions and vector bundles.} (Russian) 
Dokl. Akad. Nauk SSSR {\bf 253} (1980), no.~6, 1293--1297. 


\bibitem{ErdosV}
L.~Erd\"os, V.~Vougalter, 
{\it Pauli operator and Aharonov-Casher theorem for measure valued magnetic fields},  Comm. Math. Phys. 225 (2002), no. 2, 399--421.


\bibitem{Hille}
E. Hille, {\it Analytic function theory. Vol. II}, Ginn and Co., Boston, Mass., 1962

\bibitem{HornbergerS}
K.~Hornberger and U.~Smilansky,
{\it Magnetic edge states,}
Physics Reports \textbf{367} (2002), 249--385. 

\bibitem{Iwatsuka}
A.~Iwatsuka
{\it  The essential spectrum of two-dimensional Schr\"odinger operators with 
perturbed constant magnetic fields},
J. Math. Kyoto Univ.  \textbf{23} no. 3 (1983), 475--480.

\bibitem{Ivrii}
V.~Ivrii, {\it Microlocal analysis and precise spectral asymptotics}, Springer, Berlin, 1998


\bibitem{Landkof}
N. S. Landkof, {\it Foundations of modern potential theory}, 
Springer, New York, 1972


\bibitem{MelgardRozenblum}
M. Melgaard\ and\ G. Rozenblum, 
{\it  Eigenvalue asymptotics for weakly perturbed Dirac and Schr\"odinger 
operators with constant magnetic fields of full rank},
Comm. Partial Differential Equations {\bf 28} (2003), no.~3-4, 697--736.

\bibitem{Parfenov}
O. G. Parf\"enov, 
{\it The widths of some classes of entire functions},
Mat. Sb. {\bf 190} (1999), no. 4, 87--94; translation in Sb. Math. {\bf 190} (1999), no.~3-4, 561--568.

\bibitem{Raikov1}
G.~D.~Raikov, 
\emph{Eigenvalue asymptotics for the Schr\"odinger operator with homogeneous magnetic potential and decreasing electric potential. I. Behaviour near the essential spectrum tips,}
Comm. Partial Differential Equations {\bf 15} (1990), no.~3, 407--434;
Errata: Comm. Partial Differential Equations {\bf 18} (1993), no. 11, 1977--1979.

\bibitem{Raikov2}
G.~D.~Raikov, 
\emph{Border-line eigenvalue asymptotics for the Schr\"odinger operator with electromagnetic potential,}
Integral Equations Operator Theory {\bf 14} (1991), no. 6, 875--888.

\bibitem{Raikov3}
G.~D.~Raikov,
\emph{Spectral asymptotics for the perturbed 2D Pauli operator 
with oscillating magnetic fields. I. Non-zero mean value of the magnetic
field,}
Markov Processes Relat. Fields {\bf 9}, 775-794 (2003).


\bibitem{RaikovWarzel}
G. D. Raikov\ and\ S. Warzel, 
{\it Quasi-classical versus non-classical spectral asymptotics 
for magnetic Schr\"odinger operators with decreasing electric potentials},
Rev. Math. Phys. {\bf 14} (2002), no.~10, 1051--1072.

\bibitem{RozS}
G.~Rozenblum\ and\ N.~Shirokov, 
{\it Infiniteness of zero modes for the Pauli operator with 
singular magnetic field,}
preprint 2005, http://lanl.arxiv.org/abs/math-ph/0501059


\bibitem{Sobolev}
A.~V.~Sobolev, 
\emph{Asymptotic behavior of energy levels of a quantum particle in a homogeneous magnetic field perturbed by an attenuating electric field. I.} (Russian)  Probl. Mat. Anal., 
{\bf 9}, 67--84, Leningrad. Univ., Leningrad, 1984. English translation in:
J. Sov. Math. \emph{35} (1986), 2201--2212.

\bibitem{Sobolev2}
A.~V.~Sobolev, 
\emph{On the Lieb-Thirring estimates for the Pauli operator,}
Duke Math. J. {\bf 82} (1996), no. 3, 607--635.



\bibitem{Solnyshkin}
S.~N.~Solnyshkin, 
\emph{Asymptotic behavior of the energy of bound states of the Schr\"odinger operator in the presence of electric and homogeneous magnetic fields.} (Russian) Probl. Mat. Fiz., {\bf 10}, 266--278, Leningrad. Univ., Leningrad, 1982. 

\bibitem{StahlTotik}
H. Stahl\ and\ V. Totik, 
{\it General orthogonal polynomials}, 
Cambridge Univ. Press, Cambridge, 1992.

\bibitem{Tamura}
H.~Tamura, 
\emph{Asymptotic distribution of eigenvalues for Schr\"odinger operators with homogeneous magnetic fields,}  Osaka J. Math. {\bf 25} (1988), no. 3, 633--647.
\end{thebibliography}
\end{document}